\documentclass[10pt, a4paper, twoside]{amsart}

\usepackage{amsfonts,amssymb,amsmath,graphicx, enumitem, subfigure}
\usepackage[english]{babel}
\usepackage{booktabs}
\usepackage{multirow}
\usepackage{tikz}
\usetikzlibrary{arrows}
\usepackage{pdflscape}
\usepackage{url}
\usepackage[norelsize, ruled]{algorithm2e}
\date{\today}

\title
{Distributed Model Predictive Control for Energy Systems in Microgrids}

\author{Paul Stadler, Araz Ashouri, and Francois Mar\'echal}%
\thanks{The authors are with
         the Swiss Federal Institute of Technology Lausanne (EPFL), Sion, Switzerland.
Emails: \texttt{paul.stadler@epfl.ch, 
ashouri@alumni.ethz.ch, francois.marechal@epfl.ch}}

\begin{document}

\maketitle

\begin{abstract}
 This paper presents a flexible and modular control scheme based on distributed
 model predictive control (DMPC) to achieve optimal operation of decentralized
 energy systems in smart grids. The proposed approach is used to coordinate 
 multiple distributed energy resources (DERs) in a low voltage (LV) microgrid
 and therefore, allow virtual power plant (VPP) operation.  A sequential and 
 iterative DMPC formulation is shown which incorporates global grid targets 
 along with the local comfort requirements and performance indices. 
 The preliminary results generated by the simulation of a studied case 
 proves the benefits of applying such a control scheme to a benchmark 
 low voltage microgrid. 
\end{abstract}

\section{Introduction}
The progressive shift towards decentralized generation in power distribution 
networks has rendered the problem of optimal operation of distributed energy 
resources (DERs) to be increasingly constraining. Indeed, the integration of 
flexible deterministic energy systems combined with the strong penetration of
 uncontrollable and stochastic renewable energy sources (RESs) has pushed the
  gird to its operating limits \cite{Varaiya2011}. Micro grids represent a 
  promising concept to face the latter issue; by interconnecting and thus 
  monitoring the different DERs, the aggregated operation strategy enables 
  the provision of voltage support and other ancillary services to the local 
  distribution system operator (DSO) 
  \cite{Kennel2013, Menon2014, Bosetti2009, Katiraei2008}.
   Nevertheless, the growing integration of polygeneration systems (e.g. 
   combined heat and power units) and heat pumps raises the need of considering
    both electrical and thermal power requirements while establishing optimal
     operation strategies.

Several studies have proposed robust control mechanisms to achieve stable
 operation of low-voltage (LV) microgrids to cope with the abovementioned
  objectives. From the DSO perspective, centralized hierarchical control
   schemes are commonly selected (e.g. \cite{Bosetti2009, Hatziargyriou2005})
   , establishing set points
    for the DER load profiles to reach the grid operational targets. The 
 results presented in \cite{Borghetti2010} demonstrated the advantage of 
 implementing the following control strategy; the voltage profiles at the
  grid buses are
  highly improved through optimal DERs scheduling. However, regarding 
   the interest of DER owners connected to the microgrid, local strategies
    may highly differ from the DSO operating objective \cite{Hatziargyriou2005}.
     From the 
    point--of--view of the end-user (e.g. a residential building), model 
    predictive control (MPC) architectures are usually applied to minimize
    operating costs while satisfying the different comfort and service
     requirements of the smart grid actor 
     (e.g. \cite{Kriett2012, Collazos2009, Dagdougui2012, Henze2005}).
  Nevertheless, a lack of communication between the different MPC units while 
   defining the local strategies might cause critical grid operating
    conditions \cite{Menon2014}. In order to face these issues, authors in
     \cite{Costanzo2013} 
    presented a distributed model predictive control (DMPC) architecture,
      coordinated through an independent system operator (ISO) in
    order to steer towards a global grid objective of peak-shaving.
     However, a robust control regarding the violation of power flow
    limit in the LV network slack bus and variable bound setting 
     are not discussed in that study. 

This paper presents a flexible, modular and robust control architecture 
based on an DMPC problem formulation \cite{Richards2007}. The fully connected, 
iterative and independent MPC algorithm establishes the optimal set points of the 
 controllable DERs connected to the LV microgrid. Besides optimizing the 
 local control objectives, the coordination of the regulators allows the 
 ISO to act as a virtual power plant (VPP) while providing ancillary 
 services to the DSO. Indeed, by supplying day-ahead load predictions -- 
 computed considering a cost performance index -- to the DSO, it is able 
 to act with more certainty in the electricity market. 
The structure of the paper is the following: Section II explains the 
control scheme and Section III defines the case study considered. Section
 IV shows the advantages of using the proposed DMPC through simulation 
 results. Section V finally provides concluding comments about of the 
 proposed control method. 

\section{Control Algorithm}

In addition to avoiding a conflict of interests among the different DER 
owners and the local DSO \cite{Hatziargyriou2005}, 
the major advantage of the decentralized 
control remains in the computational effort required to solve the problem 
formulation \cite{Costanzo2013, Richards2007}. In fact, 
by splitting the large, centralized 
optimization process into multiple, smaller sub-problems, solving time 
and flexibility of the MPC is strongly improved. Nevertheless, in order 
to ensure that local actions are not threatening the global system 
stability, coordination between the different regulators is still 
required to satisfy specific coupling constraints. As presented in 
\cite{Papathanassiou2005}, several control schemes have been proposed 
to cope with the 
coordination challenges. Nevertheless, no generic design approach has 
yet been developed, leading to only case-specific, tailored DMPC problem 
formulations. 

\begin{figure}[b]
\centering
\tikzset{>=latex}
\begin{tikzpicture}
\draw  (3,2.5) rectangle (6,3.5); \draw (4.5,3) node {\textbf{ISO}};

\draw (0,-1.5) rectangle (2,-0.5); \draw (1,-1) node {\textbf{MPC$_\text{i-1}$}};
\draw (3.5,-1.5) rectangle (5.5,-0.5); \draw (4.5,-1) node {\textbf{MPC$_\text{i}$}}; 
\draw (7,-1.5) rectangle (9,-0.5); \draw (8,-1) node {\textbf{MPC$_\text{n}$}}; 

\draw [line width=1.5,->] (2,-1) -- (3.5,-1);
\draw [fill] (6.5,-1) circle [radius= 0.05];
\draw [fill] (6,-1) circle [radius= 0.05];
\draw [fill] (6.25,-1) circle [radius= 0.05];

\draw [dashed,<-] (4.125,-0.5) -- (4.125,2.5); \draw (4.125,1.625) node[left] {\textbf{\small{$\sum\limits_{j=1}^{i-1} \hat{u}_{j}^{l}(k) +$}}}; 
\draw (4.125,0.625) node[left] {\textbf{\small{$\sum\limits_{j=i+1}^{n} \hat{u}_{j}^{l-1}(k)$}}}; 
\draw [dashed,<-] (4.875,2.5) -- (4.875,-0.5); \draw (4.875,1) node[right] {\textbf{\small{$\hat{u}_{i}^{l}(k)$}}}; 

\draw (4.875,0) node[right] {\textbf{\small{(3)}}};
\draw (4.125,-1.875) node[right] {\textbf{\small{(2)}}};
\draw (3.375,0) node[right] {\textbf{\small{(1)}}};  

\draw [dashed,<-] (0.625,-0.5) -- (0.625,3.25) -- (3,3.25);
\draw [dashed,->] (1.375,-0.5) -- (1.375,2.75) -- (3,2.75);


\draw (2,-2.5) node {\small{Solving sequence}}; 
\draw (6.5,-2.5) node {\small{Information flow}}; 
\draw [dashed] (8,-2.5)--(8.5,-2.5);
\draw [line width=1.5] (3.5,-2.5)--(4,-2.5);

\draw [dashed,->] (8.375,-0.5) -- (8.375,3.25) -- (6,3.25);
\draw [dashed,<-] (7.675,-0.5) -- (7.675,2.75) -- (6,2.75);

\end{tikzpicture}
\caption{DMPC computation scheme.}
\label{fig:DMPC}
\end{figure}
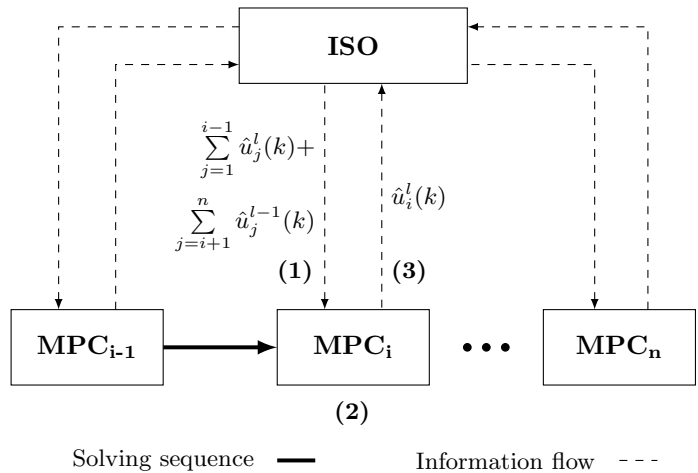

The DMPC architecture proposed in this paper relies on the combination 
of a sequential and iterative solving approach. As discussed in 
\cite{Papathanassiou2005}, regarding the information exchanged, the local 
regulators require 
comprehensive models of the different subsystems. Nevertheless, given 
the global system objective considered in this study, solely the future 
predicted control variables need to be transmitted to the different 
controllers, without any knowledge of the other DER models. Since the 
defined coupling constraint solely relates to the net power flow at 
each bus, efforts concerning the exchange of information and controller 
synchronization are strongly improved. The following paragraph details 
the specific DMPC algorithm developed for the optimal operation of DER 
units located in an LV microgrid. 

\subsection{DMPC architecture}

In the presented control algorithm, the local regulators are monitored 
through an ISO which stores the predicted load profiles of the MPCs 
connected to the microgrid. The latter independent grid agent solely 
consists of a simple data managing system, representing the aggregated 
interests of the local network actors. To solve the DMPC problem, a 
sequential updating approach \cite{Costanzo2013}
 is applied during which each regulator 
solely interacts with the ISO. As shown in Figure \ref{fig:DMPC}, the MPC hence 
calls the central ISO to gather the previously predicted operating plans 
of the remaining DERs (1), computes the local strategy (2) before 
submitting it to the ISO (3). 

This process is performed in an iterative manner until the aggregated 
load plan variance perceived by the ISO reaches the defined convergence 
criteria $\epsilon$ (Algorithm \ref{alg:DMPC}). In order to initialize the 
solving scheme for the 
first iteration (i.e. $l = 1$), the operation plan defined during the 
previous time step, $\sum \hat{u}_{}^{l}(k-1)$, is considered for units located later 
in the solving scheme. Each controller solves a standard
 MPC problem formulation. Also as shown in \cite{Jia2001}, the performance index is 
 composed of local objectives and an aggregated input target (1). The thermal 
 and electrical loads related to the local DER and service requirements 
 are computed considering a 1 hour time step (sampling time) with a 24 hours 
 prediction horizon. The control inputs $u_{i}^{}$ are defined as simply being the 
 actual state of operation/charge of the different DERs within the subsystem 
 (i.e. building) $i$. In order to prioritize thermal comfort requirements, a 
 strong penalty cost is applied when violating the soft constraint (3).

\begin{table}[b]
\begin{algorithm}[H]
\vspace{0.2cm}
\KwData{\begin{itemize} \itemsep0pt \parskip1pt \parsep0pt
\item[-] $R$ Controller set, R $\in \mathbb{N}^{*}$
\item[-] $l$ Main loop iterator  
\item[-] $U^{l}(k)$ Aggregated load profile of all controllers at time step $k$
\item[-] $\sigma_{i}^{l}(k)$ Aggregated load profile of controllers $j$, $j \in R-\{i\}$
\item[-] $\hat{u}_{i}^{l}(k)$ Predicted load profile of controller $i$ at time step $k$, $i \in R$ 
\end{itemize}
}
\KwResult{$\hat{u}_{i}(k),  i \in R$}
\textbf{Initialize:} $l = 0$, $U^{0}(k)$ = 0\;
 \While{$U^l(k) - U^{l-1}(k) < \epsilon $}{
   $l = l + 1$\;
 \For{i = 1  \KwTo n}{
        \eIf{l = 1}{
        $\sigma_{i}^{1}(k) = \sum\limits_{j=1}^{i-1} \hat{u}_{j}^{1}(k)+\sum\limits_{j=i+1}^{n} \hat{u}_{j}^{1}(k-1)$\;
         }{
        $\sigma_{i}^{l}(k) = \sum\limits_{j=1}^{i-1} \hat{u}_{j}^{l}(k)+\sum\limits_{j=i+1}^{n} \hat{u}_{j}^{l-1}(k)$\;
        }
        Compute MPC$_\text{i}(X_{i}^{l}(k),U_{i}^{l}(k),\sigma_{i}^{l}(k))$\;
        Submit $\hat{u}_{i}^{l}(k)$ to ISO\;
        }
  $U^{l}(k) = \sum\limits_{i=1}^{n}{\hat{u}_{i}^{l}(k)}$
 }
\caption{DMPC algorithm}
\label{alg:DMPC}
\end{algorithm}
\end{table}


\section{Case Study}

This case study presents the optimal operation of multiple DERs connected to 
an LV microgrid through flexible DMPC. The global objective of the distributed 
control architecture is to provide the day-ahead, aggregated load profile to 
the respective DSO of the LV network through the grid ISO. The forecasted 
consumption curve highly improves the DSO bargaining power on the daily 
electricity market since the future load profile is assumed to be determined 
\textit{a priori}. In exchange for the transmitted information, the microgrid 
end-users might benefit from a reduced electricity tariff or similar 
economic incentives from the local DSO. 

To reach the microgrid forecasted load profile, the DMPC regulators are 
however constrained to respect the aggregated consumption predictions and 
thus, need to steer towards the defined profile. In order to account for 
the stochastic related to the specific end-user behavior, the aggregated 
consumption is allowed to vary within a predefined band in which no penalty 
is perceived by the microgrid. Nevertheless, when exceeding the upper or 
lower profile limit, a non-compliance fee (Table \ref{tab:2}) is charged to the 
ISO since the economic advantage of the DSO has vanished. Hence, in the following
 case study, each grid bus (i.e. end-user) is equipped with an MPC regulator
 which optimizes the hybrid performance index (Eq. 1) which includes the local 
 operating expenses (Eq. 2) and comfort penalties (Eq. 3) in addition to the global 
 non-compliance cost (Eq. 4):  
 
\begin{equation}
\begin{cases}
J\big(X_{i}(k),U_{i}(k),\sigma_{i}(k)\big) = f_{opex}\big(U_{i}(k)\big) +
 f_{conf}\big(X_{i}(k)\big) + f_{grid}\big(U_{i}(k),\sigma_{i}(k)\big) \\  
X_{i}(k) = \{\hat{x}_{i}(k+1|k),...,\hat{x}_{i}(k+N|k)\}  \\ 
U_{i}(k) = \{\hat{u}_{i}(k|k),...,\hat{u}_{i}(k+N-1|k)\} 
\end{cases}
\end{equation}

\begin{equation}
 f_{opex}\big(U_{i}(k)\big) = 
 c_{el}(k)\big(\hat{P}_{i}^{+}(k)-\hat{P}_{i}^{-}(k)\big) +
  c_{fuel}(k)\hat{\dot{m}}_{fuel}^{+}(k)
\end{equation}
\begin{equation}
 f_{conf}\big(X_{i}(k)\big) = c_{p,conf}\big|\bar{T}_{i}(k)-\hat{T}_{i}(k)\big|
\end{equation}
\begin{equation}
 f_{grid}\big(U_{i}(k)\big) = c_{p,grid}\; \Biggl|\bar{P}_{i}^{+}(k)-
  \left( \hat{P}_{i}^{+}(k) +
   \sum\limits_{j=1,i\neq j}^{n} \hat{P}_{j}^{+}(k-1)\right)\Biggr|
\end{equation}
where the symbol {\large\^{}} indicates the predicted values of the corresponding 
state, $x$, and input, $u$, while $+/-$ superscripts indicate 
specific input and output flows (e.g. power, $P$) respectively.

\subsection{System description}
\begin{figure}[t]
    \centering
    \includegraphics[width=0.75\textwidth]{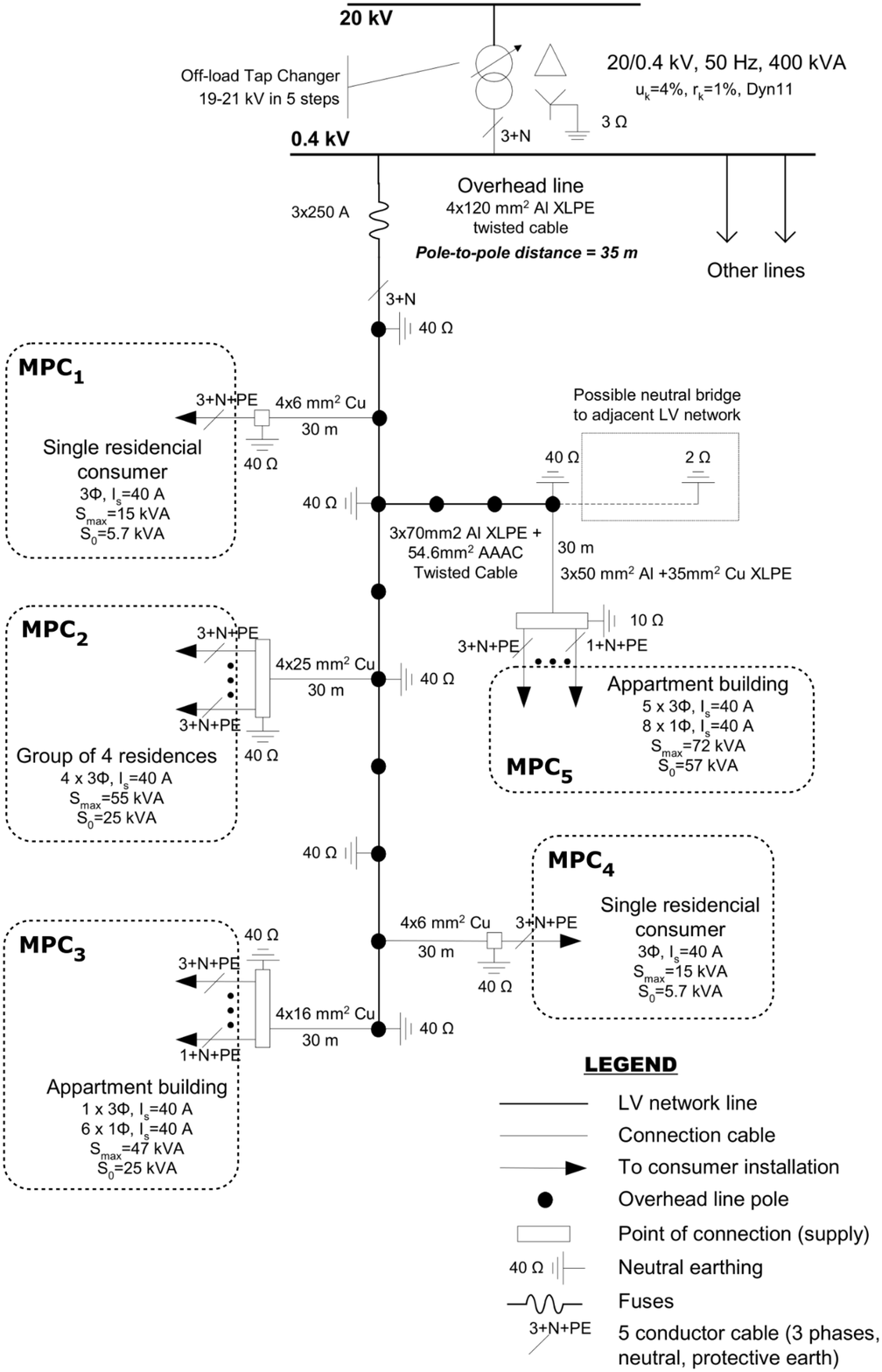}
    \caption{Low-Voltage microgrid network structure \cite{Papathanassiou2005}.}
    \label{fig:Benchmark}
\end{figure}

The system considered in this case study is composed of 
8 buildings – 6 single family houses (SFHs) and 2 
multi-family houses (MFHs) – connected to a small network. 
The different DER types located at each grid bus/buildings 
are presented in Table \ref{tab:1}. 
The simulations are performed on an IEEE benchmark LV network
 (Figure \ref{fig:Benchmark}) proposed by the CIGRE research group \cite{Papathanassiou2005}.
  The different building and DER models proposed in
  \cite{Menon2014, Collazos2009, Ashouri2014, Fazlollahi2015} have been 
  implemented in the following simulations.

\section{Results and Discussions}
This section presents the results generated by applying the proposed 
control architecture to the case study simulation introduced before.

\subsection{VPP operation -- Peak shaving} 

\begin{figure}[t]
 \centering
 \subfigure{\includegraphics[width=0.55\columnwidth]{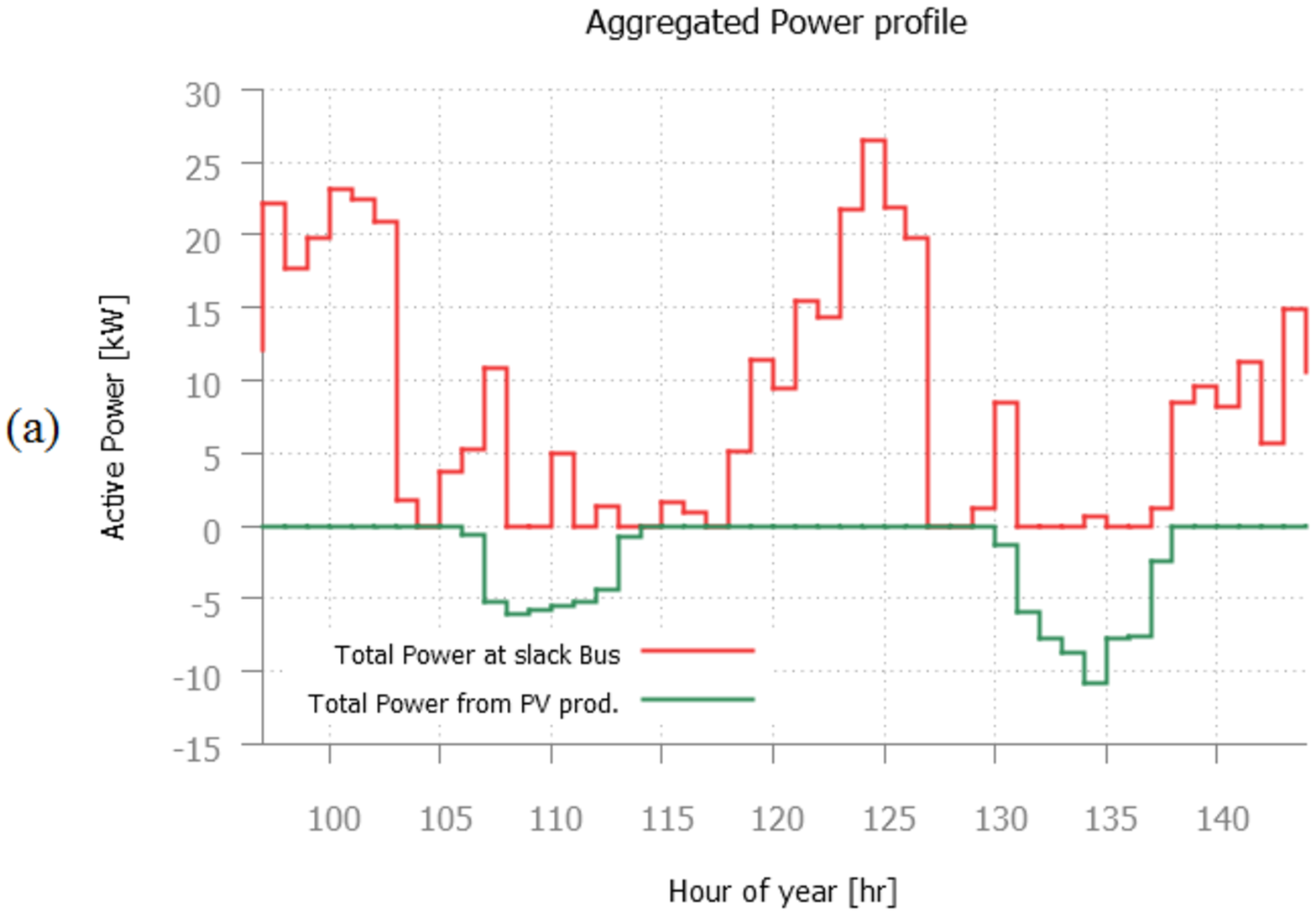}}\vspace{0.1cm} 
 \subfigure{\includegraphics[width=0.55\columnwidth]{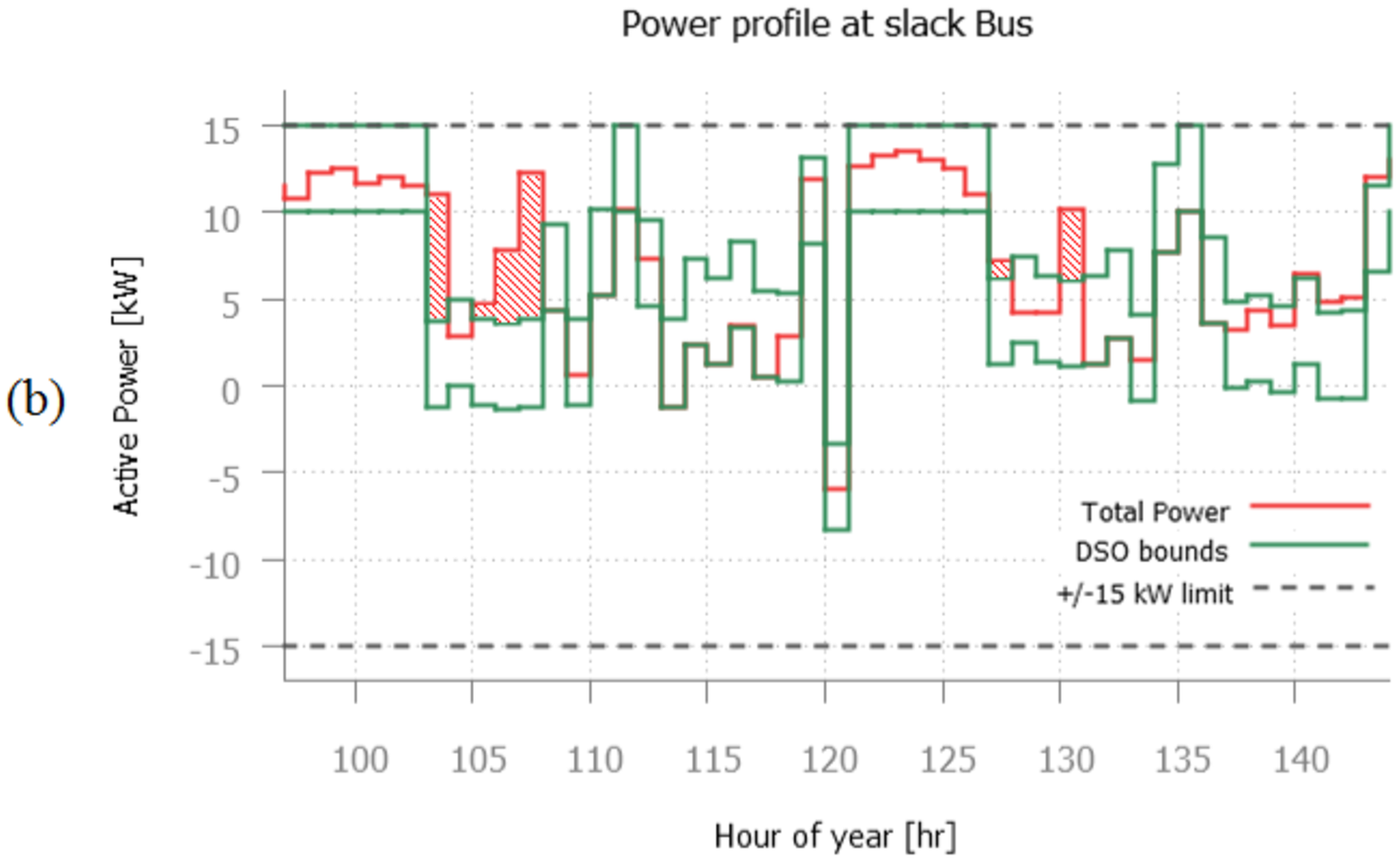}} \vspace{0cm}
 \subfigure{\includegraphics[width=0.55\columnwidth]{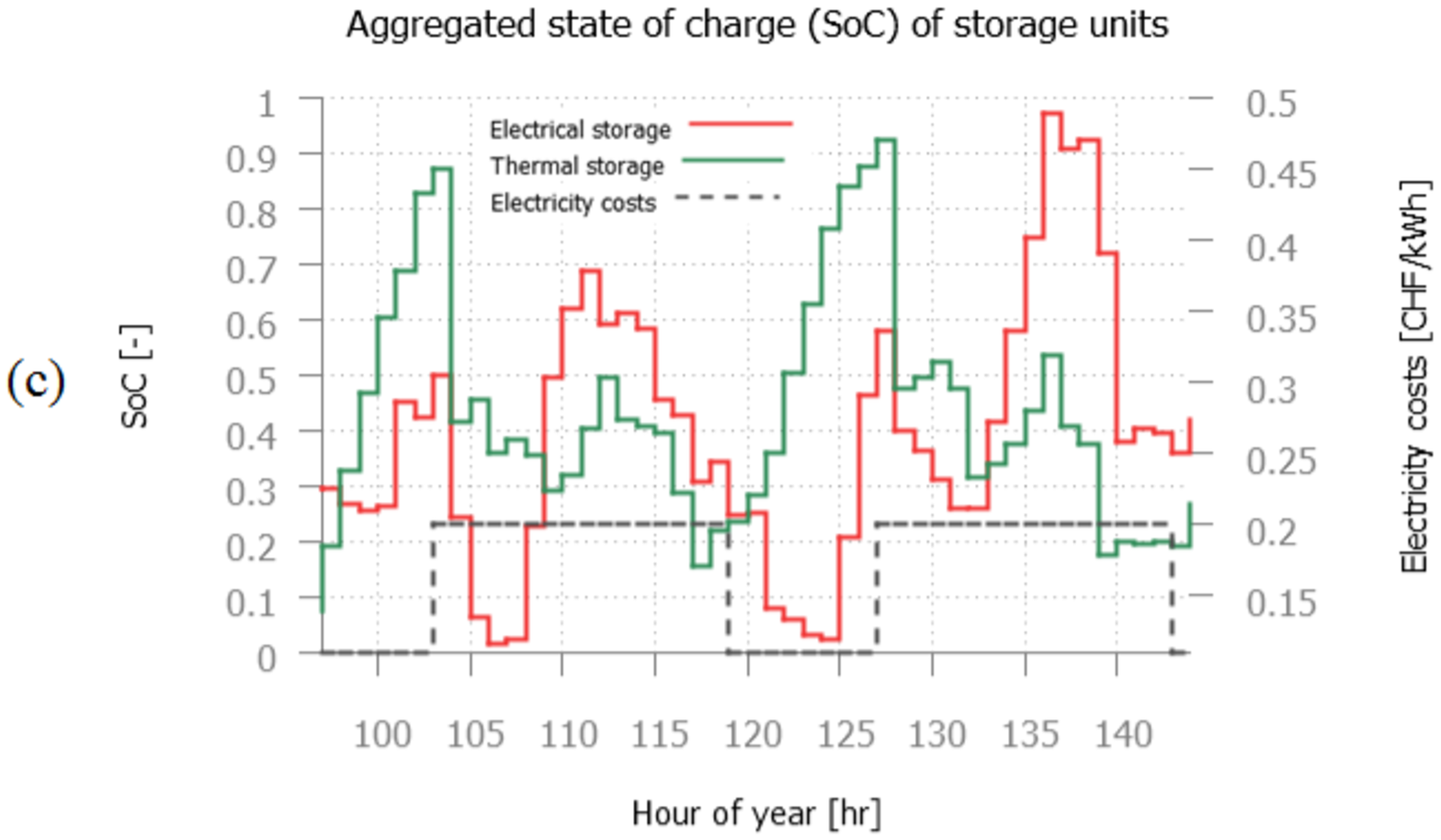}}
 \caption{MPC problem formulation comparison for 3 typical winter days.
 Active power profile at the grid slack bus without (a) and with (b) MPC 
 coordination, as well as the state of charge for aggregated storage unit (c).}
 \label{fig:win}
\end{figure}

 Figure \ref{fig:win} shows the active power flow 
 at the slack bus when the cooperation mechanism is implemented 
 (\ref{fig:win}a) and not implemented (\ref{fig:win}b) during winter time, 
 considering a 
 standard day-night electricity tariff profile (Table \ref{tab:2}). It 
 is shown that the DMPC algorithm is successfully maintaining 
 the aggregated load curve in between the predicted bounds with 
 the exception of a few overshoots. These violations are 
 correlated to strong prediction errors in the uncontrollable 
 electricity consumption which highly influence the performance 
 of MPC. Since a cost based performance index has been considered 
 (Eq. 1), the different MPC regulators tend to maximize their power 
 consumption during low electricity tariff periods (i.e. night time), 
 hence creating virtual consumption peaks.
The totally independent control MPC formulation performed without 
any coordination particularly reflects this undesired effect with 
peaks reaching $27$ kWe at the slack bus. In order to face this issue, 
a global power flow constraint of $\pm15$ kWe has been added to the 
optimization problems to attenuate the load spikes (Figure \ref{fig:win}b). 
Figure \ref{fig:win}c finally shows the aggregated state of charge (SoC) of 
the thermal and electrical storage units monitored by each MPC. 
During low tariff periods, the controllers heavily charge the 
domestic hot water (DHW) tanks and batteries to satisfy local 
(i.e. comfort and operating costs) and global (i.e. load constraints) 
objectives. Moreover, since the feed-in tariff of electricity is 
lower than the night market price (Table \ref{tab:2}), the MPC regulator 
tries to maximize self-consumption and the excess PV power generation 
is recovered by the electrical and thermal storage systems.
\begin{figure}[t]
 \centering
 \subfigure{\includegraphics[width=0.55\columnwidth]{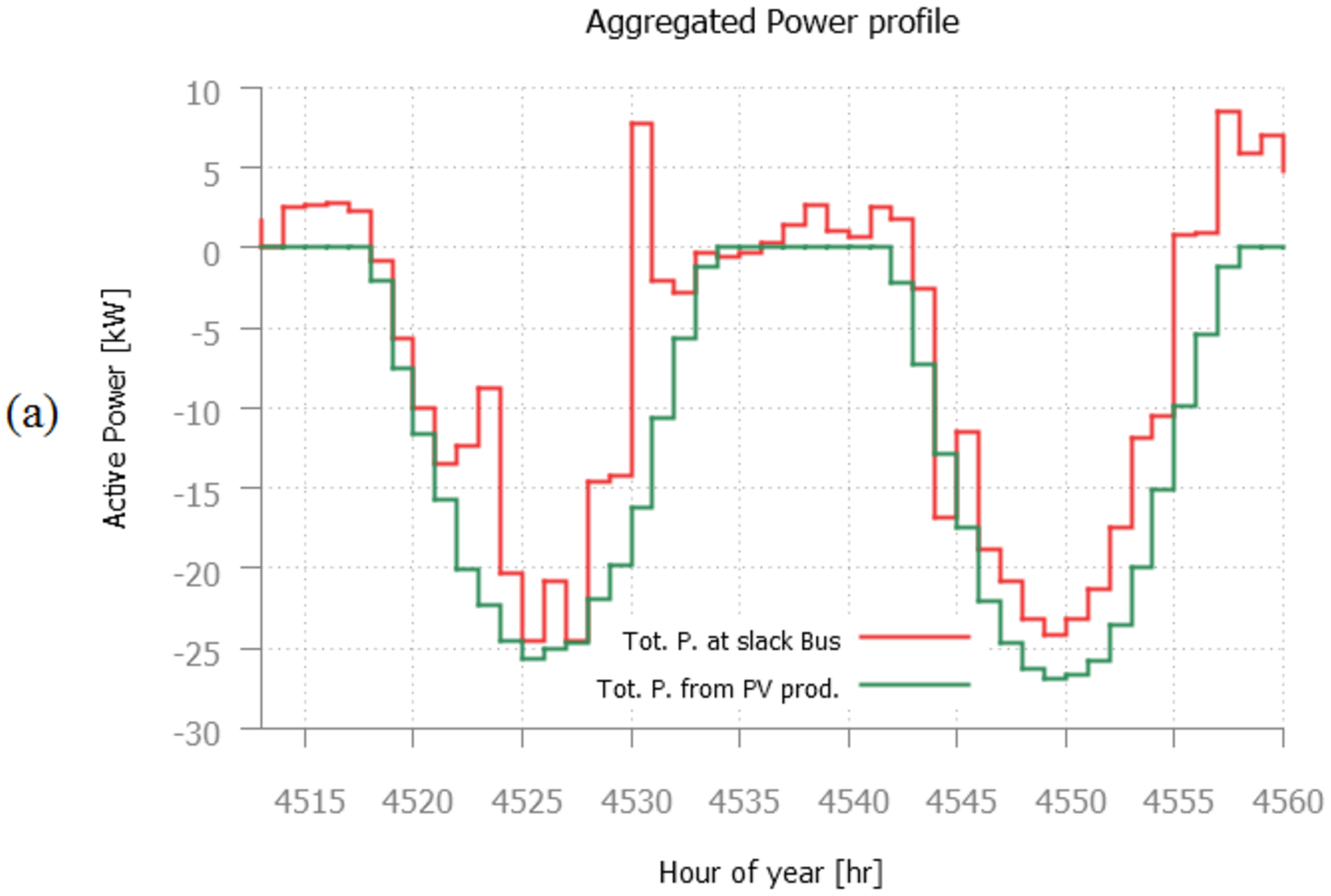}}\vspace{0.1cm} 
 \subfigure{\includegraphics[width=0.55\columnwidth]{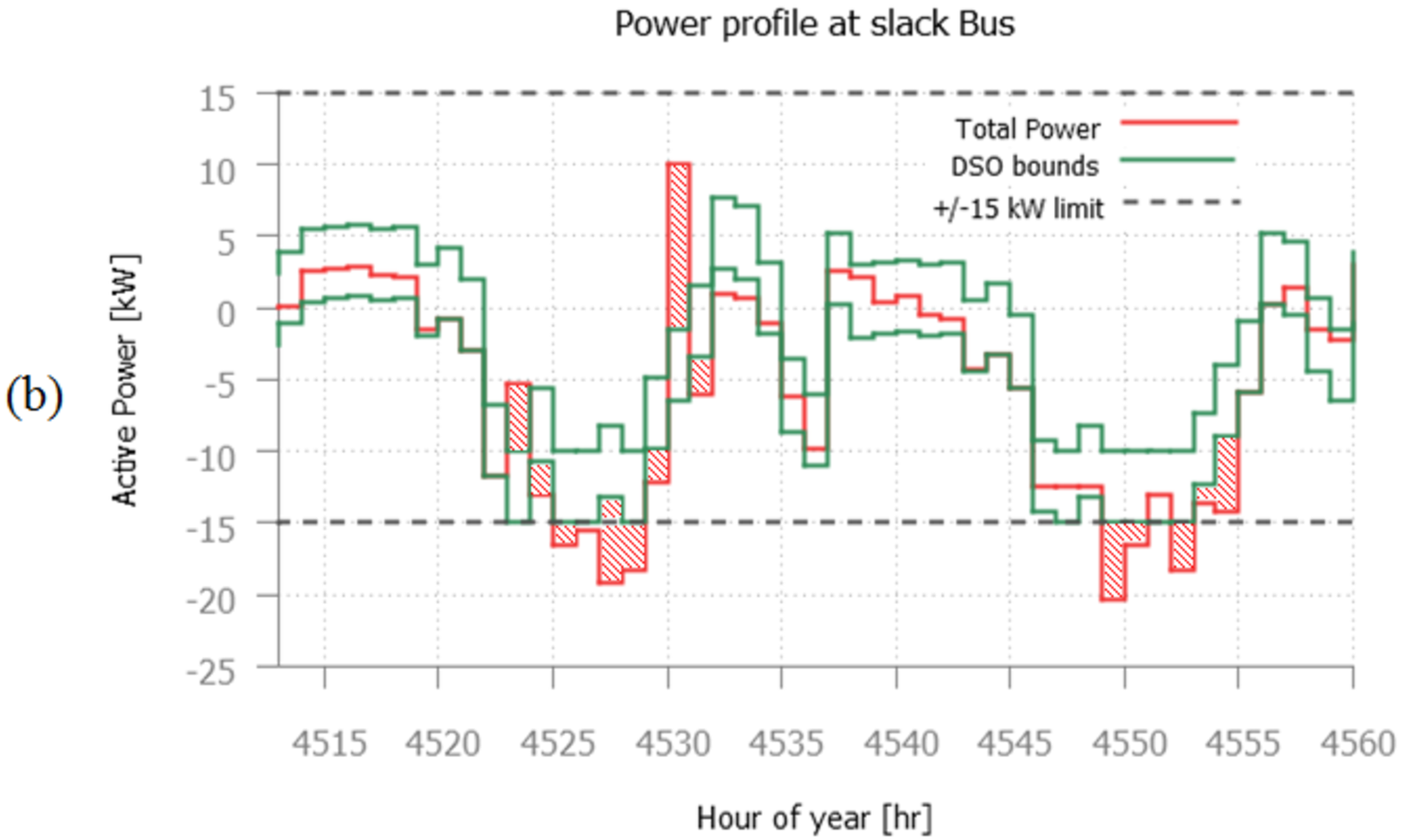}} \vspace{0cm}
 \subfigure{\includegraphics[width=0.55\columnwidth]{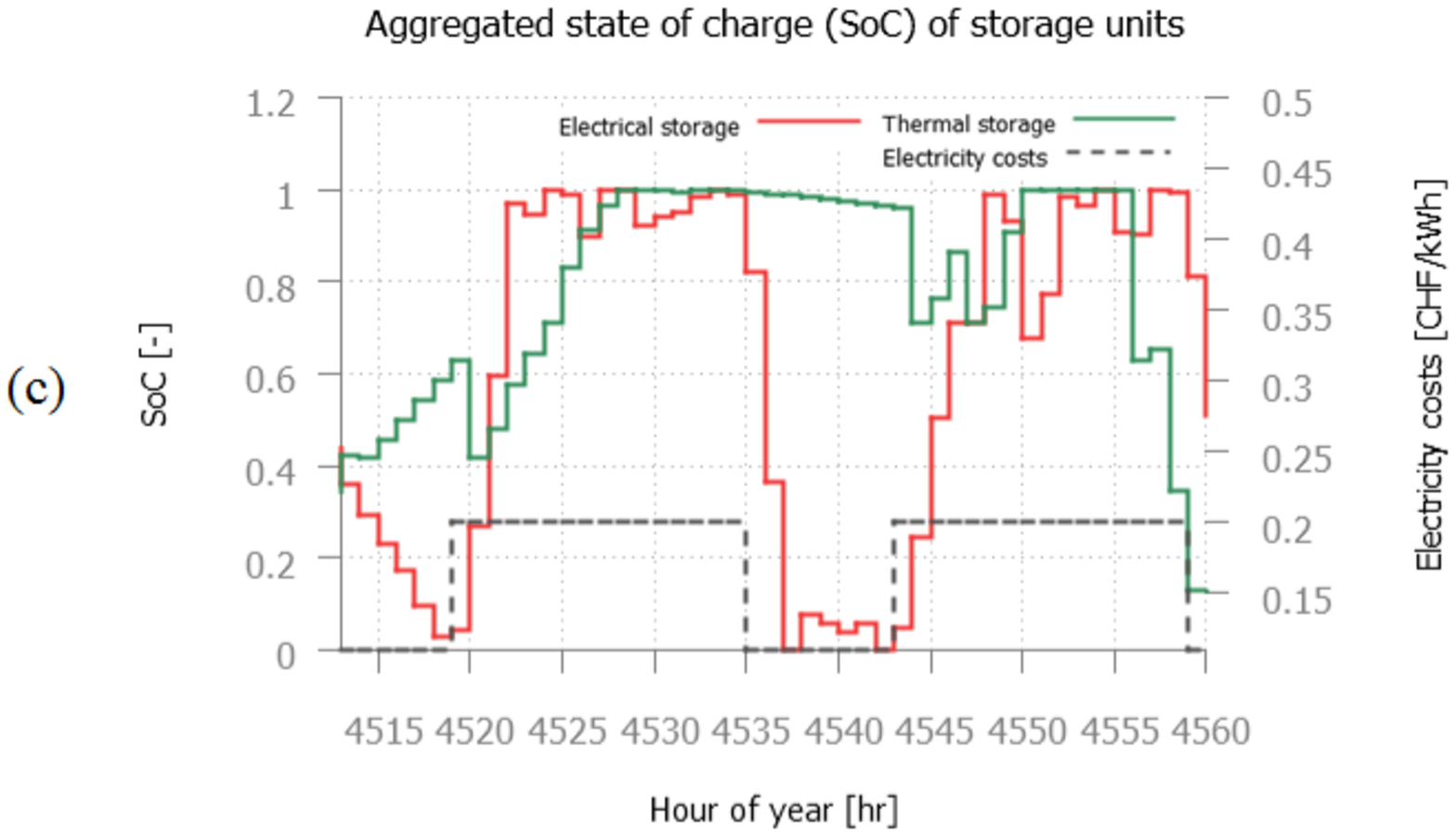}}
 \caption{MPC problem formulation comparison for 3 typical summer days.
 Active power profile at the grid slack bus without (a) and with (b) MPC 
 coordination, as well as the state of charge for aggregated storage unit (c).}
 \label{fig:sum}
\end{figure}

\begin{table}[b]
\caption{Parameters for the distributed energy system (DER).}
\label{tab:1}
\centering
\begin{tabular}{l c c c c }
\toprule
\multirow{2}{*}{\textbf{Building type}}&\multicolumn{2}{c}{\textbf{Heating}}&
\multicolumn{2}{c}{\textbf{RES}}\\
& Unit & Type & Unit & Size \\
\cmidrule(r){1-5}
SFH & Heat Pump & Air-Water & PV panel  & 50 [m$^2$]\\
MFH & Heat Pump & Air-Water & - &  \\
4 $\times$ SFH & Boiler & Natural Gas & PV panel & 100 [m$^2$]\\
SFH & Heat Pump & Air-Water & PV panel & 35 [m$^2$]\\
MFH & Cogeneration & Natural Gas & - & \\
\cmidrule(r){2-5}
&\multicolumn{2}{c}{\textbf{Thermal Storage}}&\multicolumn{2}{c}{
\textbf{Electrical Storage}}\\
& Unit & Size & Unit & Size \\
\cmidrule(r){2-5}
SFH & Hot Water Tank & 34.8 [kWh]& Battery Stack & 3 [kWh]\\
MFH & Hot Water Tank & 69.7 [kWh]& - & \\
4 $\times$ SFH & - & & Battery Stack & 12 [kWh]\\
MFH & Hot Water Tank & 46.4 [kWh]& Battery Stack & 2 [kWh]\\
SFH & Hot Water Tank & 69.7 [kWh]& - & \\
\bottomrule
\end{tabular}
\end{table}

In order to analyze the DMPC performance with respect to seasonal 
variations, Figure \ref{fig:sum} presents the active power profile when 
performing DMPC (\ref{fig:sum}b) and when considering only local objectives
 (\ref{fig:sum}a), during 
summer time. Regarding the large power generation resulting from 
the different RESs installed, the electrical storage capacity 
(Figure \ref{fig:sum}c) is not sufficient to satisfy the peak shaving constraint 
introduced previously. 
The power consumption/generation curve indeed exceeds the $\pm15$ kWe 
threshold since thermal storage through spacing heating and the hot 
water tank (that saturates over time) is not available during this 
period of the year, the excess of power must be reinjected into the 
distribution grid. To address this issue, long-term, seasonal storage 
may be necessary through a larger, centralized system located next 
to the slack bus.

\begin{table}[t]
\caption{Variable electricity tariff profiles.}
\label{tab:2}
\centering
\begin{tabular}{l c c c }
\toprule
\multirow{2}{*}{\textbf{Tariff}} &\multicolumn{2}{c}{\textbf{Import price}}&\textbf{Export price}\\
& \multicolumn{2}{c}{ \small{[CHF/kWh$_e$]}} & \small{[CHF/kWh$_e$]} \\
\cmidrule(r){1-4}
\multirow{2}{*}{Day-Night} & \textit{7 am - 10 pm} & 0.24  &\multirow{2}{*}{0.1}  \\
& \textit{else} & 0.13  & \\
\cmidrule(r){2-4}
\multirow{2}{*}{\textit{24h-ahead}} & \textit{5 pm - 11 am} & 0.21  &\multirow{2}{*}{0.1}  \\
& \textit{else} & 0.16  & \\
\cmidrule(r){1-4}
\multicolumn{4}{c}{\textbf{Violation costs} [CHF/kW$_e$]}\\
\multicolumn{2}{c}{Predicted profile} & \multicolumn{2}{c}{Global power}\\
\cmidrule(r){1-4}
\multicolumn{2}{c}{0.5} & \multicolumn{2}{c}{0.25} \\
\bottomrule
\end{tabular}
\end{table}

\subsection{24h-ahead Tariff for Load Shifting}
\begin{figure}[t]
 \centering
 \subfigure{\includegraphics[width=0.55\columnwidth]{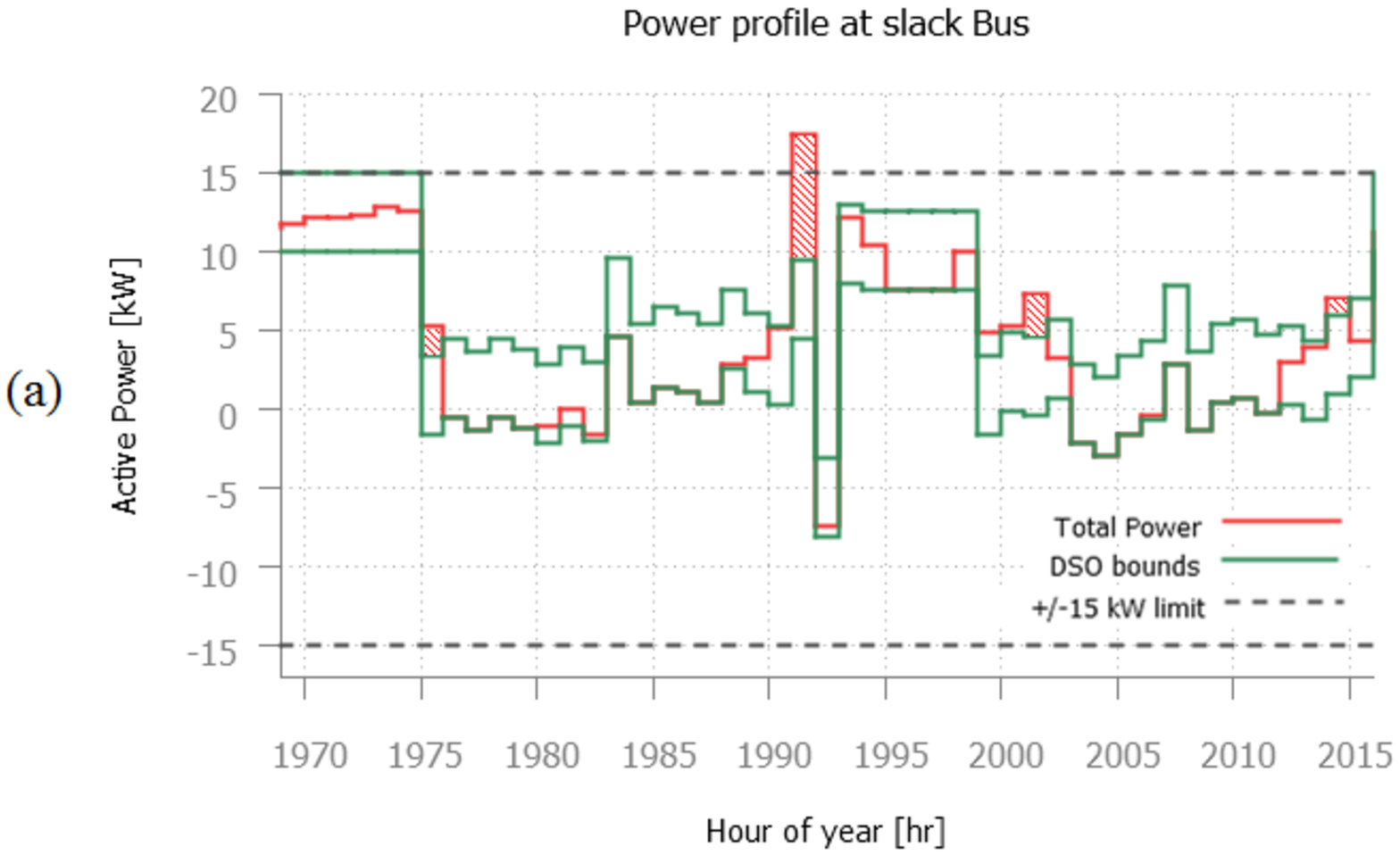}} \vspace{0.5cm}
 \subfigure{\includegraphics[width=0.55\columnwidth]{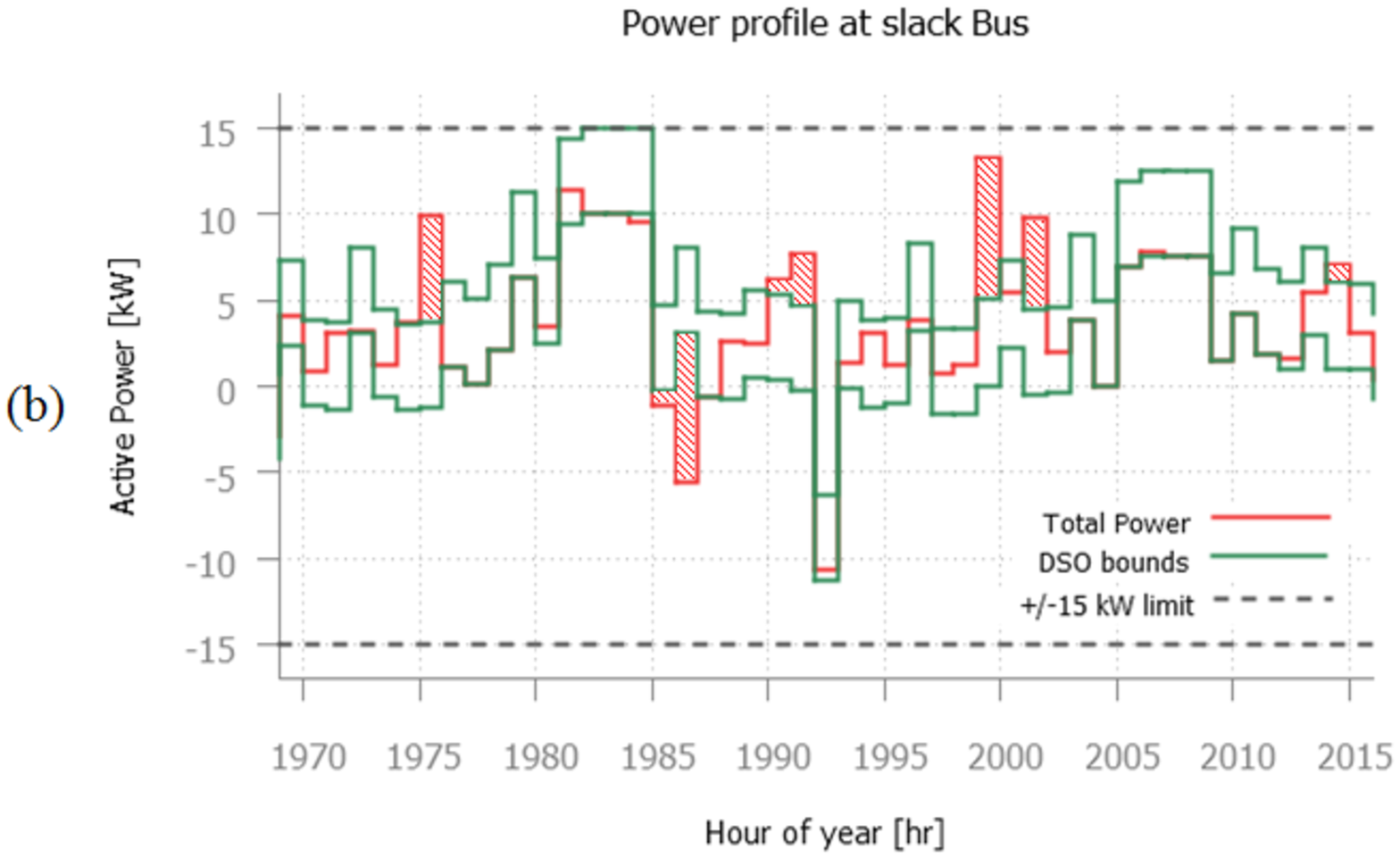}} 
 \caption{MPC problem formulation with a modified tariff profile 
 (Table \ref{tab:2}). Active power profile at the grid slack bus with the
 standard day-night (a) and the 24h--ahead (b) tariff.}
 \label{fig:ahead}
\end{figure}

As stated previously, regarding the specific performance index 
considered (Eq. 1), the different MPC controllers maximize power 
consumption during low tariff periods. Hence, in order to steer 
electricity imports towards specific periods of the day, the DSO 
might provide the ISO with specific day-ahead market prices prior 
the computation of the 24 hours load prediction. Since the MPC 
targets are entirely costs based, the local DSO may thus be able 
to shape the aggregated consumption profile regarding its current 
interest. Figure \ref{fig:ahead} shows the consumption strategy computed by the 
distributed control scheme, during spring time, while setting lower 
import costs during mid-day periods, when power generation from 
RESs is peaking (Table \ref{tab:2}). Indeed, the maximum of the aggregated 
power flow curve is shifted towards the desired time frame, thus, 
providing the DSO with an additional ancillary service, that is load 
shifting. 

Table \ref{tab:3} presents the performance values of both pricing schemes 
considered in this study (Figure \ref{fig:ahead}). As expected, the mean power 
flow (at the slack bus) during low tariff periods is higher than 
the one during low price periods. As stated previously, the total 
bound violations are related to prediction errors, mainly for the 
outer conditions. Indeed, the overshoots occurring in the morning 
time are caused by the change in the minimum comfort temperature 
within the buildings and thus, are linked to the outside temperature. 
However, the peaks observed in Figure {fig:ahead} a are considerably limited since 
the controllers used the storage systems which have been fully charged 
during the preceding low tariff period (i.e night). A similar behavior 
can be noticed for the strong undershoot in Figure \ref{fig:ahead}b; regarding 
previous low price period (i.e mid-day), the regulators could not 
provide sufficient storage capacity in order to recover the unpredicted, 
generated power excess. 

\begin{table}[b]
\caption{DMPC performance for day--night and 24h--ahead tariffs during spring time.}
\label{tab:3}
\begin{center}
\begin{tabular}{l c c c}
\toprule
\multirow{2}{*}{\textbf{Tariff}} &\multicolumn{2}{c}{\textbf{Mean power flow} \small{[kW$_e$]}} & \textbf{Total}\\
\cmidrule(r){2-3}
& \small{Low tariff} & \small{High tariff} & \small{[kWh$_e$]}    \\
\cmidrule(r){1-4}
Day-night & 10.76  &  1.78 & 202\\
\textit{24h-ahead} & 7.88  &  2.77 & 181.3\\
\cmidrule(r){1-4}
 &\multicolumn{2}{c}{\textbf{Bound violation}  \small{[kWh$_e$]}} & \textbf{Relative} \small{[-]} \\
 \cmidrule(r){2-3}
  &\small{Total} & \small{Mean} & \small{Total}  \\
\cmidrule(r){1-4}
Day-night  & 15.6 & 1.43 & 0.077\\
\textit{24h-ahead} & 34.4 & 4.82 & 0.19 \\
\bottomrule
\end{tabular}
\end{center}
\end{table}
\begin{table}[b]
\caption{Maximal voltage and angle deviation during DMPC operation.}
\label{tab:4}
\begin{center}
\begin{tabular}{l c c }
\toprule
\multirow{2}{*}{\textbf{Period}} &\multicolumn{1}{c}{\textbf{Voltage dev.}}
&\textbf{Angle dev.}\\
& \multicolumn{1}{c}{ \small{[p.u.]}} & \small{[deg]} \\
\cmidrule(r){1-3}
3-5 January & 0.01  &  2.2  \\
7-9 July & 0.007  &  4.6  \\
\bottomrule
\end{tabular}
\end{center}
\end{table}
Nevertheless, although no grid operation constraints have been 
included in the local performance, an a posteriori power flow 
analysis shows that the voltage and phase angle remain within a 
standard quality limit of $\pm 5$ p.u. (Table \ref{tab:4}). In order to account 
for the fast dynamics of electrical DERs and hence, their optimal 
operation, a second control layer should indeed be considered during 
future development. As proposed in \cite{Scattolini2009} for systems with 
different dynamics, a hierarchical approach can be added to the actual MPC 
problem formulation. The second control level evaluates the set points 
of electrical energy systems with a smaller sample time (5-15 minutes) 
in order to maintain the grid quality within the desired bounds 
($\pm 5$ p.u.). However, to cope with dynamics at higher frequencies 
(e.g. solar irradiance variations), real-time control must still 
be implemented, which relies outside of the scope of this paper. Curious reader
is referred to   
\cite{Bernstein2015}.

\section{Conclusions}
This paper has presented flexible and modular control architecture 
to define the optimal operation plans of DERs connected to a microgrid. 
Based on an iterative, independent and fully connected MPC problem 
formulation, the presented approach enables VPP operation through a 
central ISO, while providing ancillary services to the local DSO. 
The case study simulation results showed the performance of the proposed 
control scheme; although a few overshoots have been noticed particularly 
due to strong prediction errors, the coordinating MPC regulators have 
successfully steered towards the desired operation plan. In addition, 
the latter simulation showed the flexibility of the DMPC to steer the 
aggregated power profile towards specific periods through day-ahead 
market tariffs defined by the DSO. 
  
Future studies are planned to investigate the optimal design of the 
energy systems installed at each building with respect to the presented 
case study which has not been shown in this paper. Moreover, pricing 
strategies related to transgressions of the forecasted load curve might 
also be proposed. In fact, in case of constraint violation, the resulting 
penalty is preserved by the ISO, however the latter entity solely represents 
the aggregated interest of the different DERs and the DSO without providing
 any services or operating strategies.



\end{document}